\documentclass[letterpaper, reqno]{amsart}
\usepackage[T1]{fontenc}
\usepackage{amssymb, amsmath, amsthm, graphicx, mathtools, dsfont}
\usepackage{setspace} % line spacing
\usepackage[numbers]{natbib}
\usepackage[all]{xy}
\usepackage{hyperref}

\newcommand{\V}{\mathcal{V}}

\newcommand{\R}{\mathbb{R}}

\newcommand{\G}{\textsc{g}}

\newcommand{\sS}{\mathbb{S}}
\newcommand{\C}{\mathbb{C}}
\newcommand{\hH}{\mathbb{H}}

\newcommand{\Ii}{\frak{I}}

\newcommand{\inv}{^{\text{-}1}}

\newcommand{\lie}{\frak{g}}

\newcommand{\del}[1]{\frac{\partial}{\partial #1}}
\newcommand{\dev}[1]{\frac{d}{d #1}}

\newcommand{\CvxRad}{\textup{CvxRad}}
\newcommand{\HolRad}{\textup{HolRad}}

\numberwithin{equation}{section} 

\newtheorem{theorem}{Theorem}[section]
\newtheorem{lem}[theorem]{Lemma}
\newtheorem*{thma}{Theorem A}
\newtheorem*{thmb}{Theorem B}
\newtheorem*{thmc}{Theorem C}

\newtheorem{proposition}[theorem]{Proposition}
\newtheorem{cor}[theorem]{Corollary}
\newtheorem{co}[theorem]{Conjecture}

\theoremstyle{definition}
\newtheorem{df}[theorem]{Definition}
\newtheorem{exa}[theorem]{Example}

\newtheorem{pro}[theorem]{Property H}

\theoremstyle{remark}
\newtheorem{rem}[theorem]{Remark}
\newtheorem{que}{Question}

\newcommand{\bd}{\begin{df}}
\newcommand{\ed}{\end{df}}

\newcommand{\bq}{\begin{que}}
\newcommand{\eq}{\end{que}}

\newcommand{\bl}{\begin{lem}}
\newcommand{\el}{\end{lem}}

\newcommand{\br}{\begin{rem}}
\newcommand{\er}{\end{rem}}

\newcommand{\bt}{\begin{theorem}}
\newcommand{\et}{\end{theorem}}

\newcommand{\bc}{\begin{cor}}
\newcommand{\ec}{\end{cor}}

\newcommand{\bco}{\begin{co}}
\newcommand{\eco}{\end{co}}

\newcommand{\bp}{\begin{proposition}}
\newcommand{\ep}{\end{proposition}}

\newcommand{\brm}{\begin{rem}}
\newcommand{\erm}{\end{rem}}

\newcommand{\be}{\begin{equation}}
\newcommand{\ee}{\end{equation}}
\newcommand{\bx}{\begin{exa}}
\newcommand{\ex}{\end{exa}}
\newcommand{\bpr}{\begin{pro}}
\newcommand{\epr}{\end{pro}}

\begin{document}

\title{Holonomic spaces}
\author{Pedro Sol\'orzano}
\begin{abstract} A {\em holonomic space} $(V,H,L)$ is a normed vector space, $V$, a subgroup, $H$, of $Aut(V, \|\cdot\|)$ and a group-norm, $L$, with a convexity property.  We prove that with the metric $d_L(u,v)=\inf_{a\in H}\left\{\sqrt{L^2(a)+\|u-av\|^2}\right\}$, $V$ is a metric space which is locally isometric to a Euclidean ball.  Given a Sasaki-type metric on a vector bundle $E$ over a Riemannian manifold, we prove that the triplet $(E_p,Hol_p,L_p)$ is a holonomic space, where $Hol_p$ is the holonomy group and $L_p$ is the length norm defined within. The topology on $Hol_p$ given by the $L_p$ is finer than the subspace topology while still preserving many desirable properties. Using these notions, we introduce the notion of {\em holonomy radius} for a Riemannian manifold and prove it is positive. These results are applicable to the Gromov-Hausdorff convergence of Riemannian manifolds. 
\end{abstract} 

\maketitle
\section*{Introduction}
The differential geometric study of total spaces of bundles over Riemannian manifolds has been quite a prolific and influential topic. Among those bundles, the tangent bundle has been central and extensively studied by \citet{MR0112152}, \citet{MR946027}, \citet{MR2106375}, \citet*{0703059}, among many others. 

The motivation for this communication is to investigate further the properties of such bundles from the viewpoint of Metric Geometry, that is the synthetic properties of their induced metric-space structure. It is a well-known fact that the differential geometry of the tangent bundle with its Sasaki metric (or more generally any bundle with a Sasaki-type metric) is fairly rigid. In particular, its fibers (i.e. the individual tangent spaces) are totally geodesic and flat (see \cite{MR946027}). There is more than meets the eye. 

At any given point on a Riemannian manifold there are three pieces of information that interplay: the tangent space, as a normed vector space $V$; the holonomy group, as a subgroup $H$ of the isometry group of the fiber; and a group-norm $L$ on the holonomy group, given by considering the infimum 
\be\label{inlendef} 
L(a)=\inf_{\gamma}\ell(\gamma)
\ee
of the lengths of the loops $\gamma$ that yield a given holonomy element $a$. 

A {\em holonomic space} is a triplet $(V,H,L)$ consisting of a normed vector space $V$; a group $H$ of 
linear isometries of $V$; and a group-norm $L$ on such group; satisfying a local convexity property that relates them: For any element $u\in V$ there is a ball around it such that for any two elements $v,w$ in that ball the following inequality holds: 
\be
\|v-w\|^2-\|av-w\|^2\leq L^2(a)
\ee
for any element $a\in H$. See Definition \ref{holsp}.

By considering the following distance function, $d_L:V\times V\rightarrow\R$,
\be
d_L(u,v)=\inf_a\sqrt{L^2(a)+\|au-v\|^2},
\ee 
one gets a modified metric-space structure on $V$ that sheds light on the definition of a holonomic space:

\begin{thma} A triplet $(V,H,L)$  is a holonomic space if and only if $d_L$ is locally isometric to usual distance on $V$.
\end{thma}

This is Theorem  \ref{lociso}, proved in section \ref{Hols}. 

The measure of nontriviality of a holonomic space is controlled by the {\em holonomy radius}, a continuous function on $V$, given by the supremum of the radii of balls for which the local convexity property is satisfied. This function is finite if and only if $H$ is nontrivial [Proposition \ref{flatinfty}].

Considering the holonomy radius at the origin already yields some information on the group-norm in the case when the normed vector space is actually an inner product space. Namely the following result.

 \begin{thmb} Given a holonomic space $(V,H,L)$, the identity map on $H$ is Lipschitz between the left invariant metrics on $H$ induced by $L(a\inv b)$ and $\sqrt{2\|a-b\|}$ respectively, where $\|\cdot\|$ stands for the operator norm. Moreover, the dilation is precisely the reciprocal of the holonomy radius $\rho_0$ at the origin of $V$.
\be
 \sqrt{2\|a-b\|}\leq\frac{1}{\rho_0}L(a\inv b).
\ee
\end{thmb}
 This is a consequence of Theorem \ref{holrad0} and Corollary \ref{lipnorm} in section \ref{Hols}.

Recall that a Sasaki-type metric $\G$ on a Euclidean vector bundle with compatible connections is given in terms of the connection map $\kappa:TE\rightarrow E$, uniquely determined by requiring that $\kappa(\sigma_*x)=\nabla^E_x\sigma$, as
\be
\G(\xi,\eta)=g(\pi_*\xi,\pi_*\eta)+h(\kappa\xi,\kappa\eta),
\ee
for vectors $\xi,\eta\in TE$.

Given these considerations one gets the following result.
\begin{thmc} Given a Euclidean vector bundle with a compatible connection over a Riemannian manifold, each point in the base space has a naturally associated holonomic space, with the fiber over that point being the underlying normed vector space.

Furthermore, if the total space is endowed with the corresponding Sasaki-type metric then the aforementioned modified metric-space structure coincides with the restricted metric on the fibers from the metric on the total space.  
\end{thmc}
This result is stated more precisely in Proposition \ref{distsas}, Theorem \ref{lengthproof}, and Theorem \ref{holfib}.

The group-norm in Theorem C is given by \eqref{inlendef}. The study of this group-norm was already hidden in the work of \citet{MR1793685} and \citet{MR1125864}.

This group-norm induces a new topological group structure on the holonomy group that makes the the group-norm continuous while retaining the continuity of the holonomy action [Lemma \ref{conthol}]. In should be noted that with the standard topology (i.e. that of a Lie group) of the holonomy group, this group-norm is not even upper semicontinuous. \citet{MR1125864} had already noted this (an immediate example is to consider a metric that is flat in a neighborhood of a point and consider the group-norm associated at that point). He proved that if the Lie group topology is compact then the group-norm topology is bounded, which is a surprising result given that the group-norm topology is finer. 

\citet{MR1793685} defines a `size' for a given holonomy element as an infimum over {\em acceptable} smooth metrics on the holonomy group (quoted here as \ref{Tappsize}). As such, he proved that holonomy `size'  and the length group-norm \eqref{inlendef} are comparable up to a constant that depends only on the base space and the norm of the curvature (see Theorem \ref{Tapphol}). These results are discussed in more detail in the last section of this communication.  

Finally, as an application of this constructions, the {\em holonomy radius} of a Riemannian manifold is defined by assigning to each point the holonomy radius at the origin of the holonomic space associated to it by Theorem C in the case of its tangent bundle with the Levi-Civita connection [Definition \ref{holradmfl}].  

In view of Theorem B, this function is positive [Theorem \ref{posradmfl}] and has a precise formulation in terms of the group-norms and of the usual linear action of the holonomy groups. 

%%%%%%%%%%%%%%%%%%%%%%%%%%%%%%%%%%%%%%%%
\section{Background on Normed Groups and Semimetrics}\label{Top} 
The definitions and results stated in this section are here for completeness and for the sake of setting the background and notation for the sequel. The notion of semimetrics (often called pseudometrics in the literature) is reviewed, so is the concept of (left) isometric group actions. In this generality, quotient spaces are not necessarily well behaved, but further quotients can be taken to return to the category of metric spaces. These quotients will however be related. 

Since the concepts needed in following sections are mainly metric, the topological and smooth structural results will be omitted or freely used whenever needed.  
\subsection{Norms and distances}
\bd\label{grpnm} Let $G$ be any group. A {\em group-norm} on $G$ is a function $N:G\rightarrow \R$ that satisfies the following properties.
\begin{enumerate}
\item Positivity: $N(A)\geq0$
\item Non-degeneracy: $N(A)=0$ iff $A=id_V$
\item Symmetry: $N(A\inv)=N(A)$
\item Subadditivity (``Triangle inequality''): $N(AB)\leq N(A)+N(B)$.
\end{enumerate}
\ed
\bx\label{nrmex1} Let $V$ be a normed vector space and let $G$ be a subgroup of the  group of norm preserving automorphisms of $V$. Then $N(A)=\|id_V-A\|$ is a group-norm.
\ex
\bx\label{nrmex2} Let $f:[0,\infty)\rightarrow [0,\infty)$ be any non-decreasing subadditive function, $f(t+s)\leq f(t)+f(s)$, with $f(0)=0$. Let $N:G\rightarrow\R$ be any group-norm on $G$. Then $f\circ N$ is also a group-norm on G.
\ex
\bp\label{grpnrmtop} A group $G$ together with a group-norm $N$ becomes a topological group with the left invariant metric induced by 
\be d(A,B)=N(A\inv B).
\ee
\ep
\begin{proof} Left-invariance follows from the fact that $(CA)\inv CB=A\inv B$.  The map $(A,B)\mapsto A\inv B$ is continuous since
\begin{eqnarray*}
d(A\inv B,C\inv D)&=&N(B\inv AC\inv D)\\
&\leq& N(A\inv B)+N(C\inv D)=d(A,B)+d(C,D)\\
&\leq&\sqrt{2}\sqrt{d^2(A,B)+d^2(C,D)} .
\end{eqnarray*}
\end{proof}

\bd  Given a group-norm $N$ on a group $G$, the topology generated by $N$ will be called {\em the $N$-topology} on $G$.  
\ed

\bp\label{contnorm} With the $N$-topology on $G$, the group-norm $N$ is continuous. 
\ep
\begin{proof} This follows from the fact that $|N(A)-N(B)|\leq N(A\inv B)$, which in turn follows directly from the triangle inequality in Definition \ref{grpnm}.
\end{proof}

As seen, the notion of group-norm is completely equivalent to that of a left-invariant metric on a group; given that it appears in examples is it considered separately. It is also equivalent to right-invariant metrics. For more details on normed groups see the survey by \citet{NORMED}.

\subsection{Semimetrics}

\bd Let $X$ be a set together with a function from $d:X\times X\rightarrow\R$. $(X,d)$ is a {\em semimetric space} if $d$ is nonnegative, symmetric and satisfies the triangle inequality. 
\ed

\bp[see \cite{MR2048350}]\label{setome} Given a semimetric space $(X,d)$, let $x\sim y$ if $d(x,y)=0$. Then $X'=X/\sim$ is a metric space with metric, $d'$,
\be\label{dblquomet}
d'([x],[y])=d(x,y)
\ee
for any choice of representatives. Also, the canonical projection map is open and continuous with the quotient topology. 
\ep

\subsection{Isometric Transformation Groups} 

\bd
An {\em isometric group action} consists of a triplet $(G,X,\varphi)$, where $G$ is an abstract group, $(X,d)$ a metric space and $\varphi: G\rightarrow \textup{Iso}(X,d)$ a group homomorphism. The {\em orbit} of a point $x\in X$, denoted by $G(x)$, is the equivalence class of all $y\in X$ such that $y=gx:=\varphi(g)(x)$ for some $g\in G$. The space of equivalence classes is called {\em orbit space}  and will be denoted by $G\backslash X$. If, furthermore, $G$ is a topological group and the map $(g,x)\mapsto gx$ is continuous then $(X,G)$ is a {\em transformation group}.
\ed

\br The quotient map from $X$ to $G\backslash X$ is an open continuous map with the quotient topology.
\er

The following fact will be used in the sequel. It is a classical result of spaces of continuous maps with the compact open topology that the continuity of the evaluation map is equivalent, under some assumptions, to the continuity of the embedding (see \citet{MR0464128}). 
 
\bp\label{COcont} If $X$ is locally compact then a transformation group $(X,G)$ is equivalent to a continuous homomorphism $\varphi: G\rightarrow \textup{Iso}(X,d)$, where the codomain has the compact open topology.
\ep

\bl[see \cite{MR1393940}. p. 47]\label{techlem} Let $(X,d)$ be a locally compact metric space and let $\{\varphi_i\}$ be a sequence of isometries of $(X,d)$. If there exists a point $x\in X$ such that $\{\varphi_i(x)\}$ converges, then there exists a subsequence $\{\varphi_{i_k}\}$ that converges to an isometry of $(X,d)$.
\el
With the hypotheses of the previous lemma, one has the following fact, which for future reference is included here.
\bp Let $(X,G)$ be a transformation group where $X$ is locally compact and connected. Then $G\backslash X$ is a semimetric space with
\be\label{grpquomet}
d(G(x),G(y))=\inf_gd(x,gy)
\ee
Furthermore, let $H$ be the closure of $\varphi(G)$ in the isometry group of $X$. Then there exists an isometry such that
\be H\backslash X\cong (G\backslash X)/\sim.
\ee
\ep
\begin{proof}
All the properties of a semimetric are straightforward. Let $x\in X$, then its equivalence class $[x]\subseteq X$ in the right-hand side of the equation is the set  
$$
[x]=\{y:\exists\varepsilon>0,\exists g\in G,d(x,gy)<\varepsilon\}.
$$
This is equivalent to $H(y)$  for any fixed $y\in[x]$ since one can produce a sequence $\{g_i\}$ of isometries such that the sequence $\{g_i(y)\}$ converges to $x$; thus  by Lemma \ref{techlem} there is a $g\in H$ with $x=gy$. There is therefore a canonical bijection between both sides of the equation. It follows that it is an isometry since the metric on each side is defined to be the distance between equivalence classes as subsets of $X$ (cf. \eqref{dblquomet}, \eqref{grpquomet}). 
\end{proof}

\bc Let $(X,G)$ be a transformation group where $X$ is locally compact and connected. The orbit space is a metric space if $\varphi(G)$ is a closed subgroup of $\textup{Iso}(X,d)$.
\ec

\bp\label{twisted} Let $(G,d_G)$ be a metric topological group with left-invariant metric $d_G$. Let $(X,G)$ be a transformation group (no assumption on connectedness or local compactness). Then so is $(X\times G, G)$, where the action is given by $g(x,h):=(gx,gh)$. The quotient space $G\backslash(X\times G)$ is a metric space and is homeomorphic to $X$ under the map $x\mapsto G(x,e)$, the orbit of $(x,e)$. This induces a new metric on $X$ given by 
$$
d'(x,y)=\inf_g\sqrt{d_G^2(e,g)+d^2(x,gy)}.
$$
\ep
\begin{proof}
That $(X\times G, G)$ is a transformation group is immediate from the hypotheses. The second quotient, $G\backslash (X\times G)/\sim$, is a metric space by Prop \ref{setome}. Let $[x,g]:=[G(x,g)]$ stand for the equivalence class  $(x,g)$ in the second quotient, regarded as a subset of $X\times G$. The map 
\be
x\mapsto[x,e]
\ee
is a injective; indeed, suppose that $[x,e]=[y,e]$. For all $\eta>0$ there exists a group element $g$ such that $d(x,gy),\ d_G(e,g)<\eta$. Then, by continuity of the action, for every $\varepsilon>0$ there exists $0<\delta\leq\varepsilon/2$ such that for all $g$, $d_G(g,e)<\delta$, $d(y,gy)<\varepsilon/2$ (by equation \eqref{grpquomet}). So by letting $\eta=\delta$,
$$
d(x,y)\leq d(x,gy)+d(y,gy)<\varepsilon.
$$
It is also onto since $[x,g]=[g\inv x,e]$. It is continuous since it is a composition of continuous maps and it is open since the quotient maps are open and the map $(x,g)\mapsto g\inv x$ is continuous, hence a homeomorphism.

Since the map $x\mapsto G(x,e)$ is also bijective and continuous, it follows that $G\backslash X$ was already a metric space, as claimed.
\end{proof}

For a systematic introduction to (smooth) tranformation groups, see the lecture notes by \citet{MR2459325}.
%%%%%%%%%%%%%%%%%%%%%%%%%%%%%%%%%%%%%%%%
\section{Holonomic spaces}\label{Hols} 
The notion of  a {\em holonomic space} is introduced in this section. It will be seen in the sequel how these spaces occur as fibers of Euclidean vector bundles with suitable conditions imposed.  Several properties of holonomic spaces are also analyzed here.
 
\bd\label{holsp} Let $(V,\|\|)$ be a normed vector space, $H\leq Aut(V)$ a subgroup of norm preserving linear isomorphisms, and $L:H\rightarrow \R$ a group-norm on $H$. The triplet $(V,H,L)$ will be called a {\em holonomic space} if it further satisfies the following convexity property: 
\begin{itemize}
\item[(P)] For all $u\in V$ there exists $r=r_u>0$ such that for all $v,w\in V$ with $\|v-u\|<r$, $\|w-u\|<r$, and for all $A\in H$,
\be\label{holspeq}
\|v-w\|^2-\|v-Aw\|^2\leq L^2(A).
\ee
\end{itemize}
\ed
\bd Let $(V,H,L)$ be a holonomic space. The {\em holonomy radius} of a point $u\in V$ is the supremum of the radii $r>$ satisfying the convexity property (P) given by \eqref{holspeq}. It will be denoted by $\HolRad(u)$. It may be infinite.
\ed

\bl\label{contid} Given a holonomic space $(V,H,L)$ as above, there exists $r>0$ such that for $u\in V$, $|u|<r$, and for any $B\in H$,
\be\label{ctid}
\|u-Bu\| \leq L(B).
\ee
\el
\begin{proof} Simply choose $r=r_0$ as in \ref{holsp}, $v=Bu$ and $A=B\inv$.
\end{proof}
\bd Given a holonomic space $(V,H,L)$, the largest radius of a ball satisfying Lemma \ref{contid} is the {\em convexity radius of a holonomic space}. It may be infinite. 
\ed
\br The convexity radius is in general larger than the holonomy radius than the holonomy radius at the origin, as can be seen in Example \ref{counterex}
\er
Recall that the group norm $L$ on $H$ induces a topological group structure on $H$, the $L$-topology (see Proposition \ref{grpnrmtop}).
\bl\label{conthol} Given a holonomic space $(V,H,L)$, the action $H\times V\rightarrow V$ is continuous with respect to the $L$-topology on $H$. Furthermore, it is uniformly continuous when restricted to bounded invariant sets.
\el 
\begin{proof}  Let $r_0$ be the convexity radius.  Let $(a,u)\in H\times V$, with $\|u\|\geq r_0$. Fix $\lambda>0$ such that $\lambda\|u\|\l< r$,  and let $\varepsilon>0$. Let $K=\sqrt{1+\frac{1}{\lambda^2}}$. Note that for any positive real numbers $x,y\in\R$,
\[x +\frac{1}{\lambda}y\leq K\sqrt{x^2+y^2}.\] 

Now, if $\delta=\min\{\frac{\varepsilon}{K},\frac{r-\lambda\|u\|}{\lambda}\}$ and $\sqrt{L^2(a\inv b)+\|u-v\|^2}<\delta$. Notice that  \[
\|\lambda v\|\leq\lambda\|u-v\|+\|\lambda u\|< \lambda\delta + \lambda\|u\| \leq r_0.
\]

Thus, 
\begin{align*}
\|au-bv\|&=\|u-a\inv bv\|=\frac{1}{\lambda}\|\lambda u - a\inv b \lambda v\| \\
&\leq\frac{1}{\lambda}\left( \|\lambda u-\lambda v\|+\|\lambda v-a\inv b\lambda v\|\right) \\
&\leq\|u-v\|+\frac{1}{\lambda}L(a\inv b)\\
&\leq K\sqrt{L^2(a\inv b)+\|u-v\|^2}\\
&= K\sqrt{L^2(a\inv b)+\|u-v\|^2}<K\delta\leq\varepsilon,
\end{align*}
\end{proof}

\bt Let $(V,H,L)$ be a holonomic space.  
\be\label{hlmet} d_L(u,v)=\inf_{a\in H}\left\{\sqrt{L^2(a)+\|u-av\|^2}\right\},
\ee
is a metric on $V$. 
\et
\begin{proof} By Lemma \ref{conthol} one sees that the action $H\times V\rightarrow V$ is continuous with respect to the $L$-topology on $H$. Letting $G=H$ and $X=V$  in Proposition \ref{twisted} it follows that the  $V$ is homeomorphic to $H\backslash(V\times H)$ and that the pullback metric on $V$ is given by \eqref{hlmet}.
\end{proof}

\bd
Given a holonomic space $(V,H,L)$. The metric given by \eqref{hlmet} will be called {\em associated holonomic metric} and $V$ together with this metric will be denoted by $V_L$.
\ed

\bt\label{lociso}A triplet $(V,H,L)$ is a holonomic space if and only if $id:V\rightarrow V_L$ is a locally isometry. 
\et
\begin{proof}
By property (P), given any point $u\in V$ there exists a radius $r>0$ such that for all $v,w\in V$, with $\|v-u\|<r$ and $\|w-u\|<r$, and for all $A\in H$, 
$$\|v-w\|\leq\sqrt{L^2(A)+\|v-Aw\|^2}.$$

 Hence, considering the infimum of the right-hand side, it follows that
$$ \|v-w\|\leq d_L(v,w)\leq\sqrt{L^2(id_V)+\|v-w\|^2}=\|v-w\|.$$ 
Conversely, if the identity is a local isometry, property (P) in Definition \ref{holsp} is also satisfied: Let $B$ be a ball around $u\in V$ on which the identity map $id_V|B$ is an isometry. Therefore, for any $A\in H$ and any pair of points $v,w\in B$, 
$$
\|v-w\|=\inf_{a\in H}\left\{\sqrt{L^2(a)+\|v-aw\|^2}\right\}\leq \sqrt{L^2(A)+\|v-Aw\|^2}.
$$
\end{proof}
\br The holonomy radius is also the radius of the largest ball so that the restricted $d_L$-metric is Euclidean.
\er
\bp Let $(V,H,L)$ be a holonomic space. The original norm on $V$ is recovered by the equation
\be
\|v\|=d_L(v,0)
\ee
\ep
\begin{proof} Because $H$ acts by isometries on $V$, 
$$
d_L(v,0)=\inf_{a\in H}\left\{\sqrt{L^2(a)+\|v\|^2}\right\}.
$$
The conclusion now follows by letting $a=id_V$.
\end{proof}
\bc Given a holonomic space $(V,H,L)$ the rays emanating from the origin are geodesics in $V_L$.
\ec
\bp\label{flatinfty} Let $(V,H,L)$ be a holonomic space. Then $H=\{id_V\}$ if and only if there exists $u\in V$ for which the holonomy radius is not finite.
\ep
\begin{proof} If $H$ is trivial, then $L\equiv0$, and so $V$ is globally isometric to $V$, hence for any $u\in V$ the holonomy radius is infinite. Conversely, if there exists $u\in V$ with $\HolRad(u)=\infty$, and there is $a\in H$ with $L(a)>0$, then for any $v\in V$, $\|v-av\|\leq L(a)$ should hold. This is a contradiction since $v\mapsto\|v-av\|$ is clearly not bounded.
\end{proof}

\bc\label{posrad} Let $(V,H,L)$ be a holonomic space. Then the function $u\mapsto \HolRad(u)$ is positive. Furthermore, it is finite provided $H$ is nontrivial.
\ec
\bp Let $(V,H,L)$ be a holonomic space. The function  $u\mapsto \HolRad(u)$ is continuous. 
\ep
\begin{proof} By Proposition \ref{flatinfty} one can assume, with no loss of generality, that $H\neq\{id_V\}$. Let $u\in V$ and let $\varrho(u)$ be the holonomy radius at $u$. Let $v\in V$ with $\|v-u\|<\varrho(u)$, i.e. $v\in B_{\varrho(u)}(u)$,  then by maximality of $\varrho(v)$, it has to be at least as large as the radius of the largest ball around $v$ completely contained in $B_{\varrho(u)}(u)$,  
$$
\varrho(v)\geq \varrho(u)-\|u-v\|.
$$ 
Also, by maximality of $\varrho(u)$, if follows that $\varrho(v)$ cannot be strictly larger than the smallest ball around $v$ that contains $B_{\varrho(u)}(u)$, 
$$\varrho(v)\leq \varrho(u)+\|u-v\|.$$ 
Therefore, at any given point $u\in V$ and any $\varepsilon>0$, there exists $\delta=\min\{\varrho(u),\varepsilon\}$ such that for any $v\in V$, $\|u-v\|\leq\delta$ it follows that 
$$
|\varrho(u)-\varrho(v)|\leq\|u-v\|\leq\varepsilon.
$$
\end{proof}
\bt\label{cvxrad0} Let $(V,H,L)$ be a holonomic space. Then the convexity radius is given by
\be\label{cvxrad}
\CvxRad=\inf_{a\in H}\frac{L(a)}{\|id_V-a\|},
\ee
where for any $T:V\rightarrow V$,  $\|T\|$ denotes its operator norm.
\et
\begin{proof} Let $u\in V$ with $\|u\|\leq\frac{L(A)}{\|id_V-A\|}$, then
$$
\|Au-u\|\leq\|A-id_v\|\|u\|\leq L(A)
$$
which proves that 
$$
\CvxRad\geq\inf_{a\in H}\frac{L(a)}{\|id_V-a\|}.
$$
Now, let $\varrho>\frac{L(A)}{\|id_V-A\|}$ and let $\varepsilon>0$ be such that
$$
\varepsilon < \|A-id_V\|\varrho-L(A)=\|A-id_V\|\left(\varrho-\frac{L(A)}{\|A-id_V\|}\right)>0.
$$
Then, by the definition of operator norm, there exists $u\in V$ with $\|u\|=\varrho$ such that
$$
\|A-id_V\|\varrho\geq \|Au-u\|>\|A-id_V\|\varrho-\varepsilon.
$$ 
The second inequality yields that
$$
\|A-id_V\|\varrho-\varepsilon>L(A).
$$
This proves that $\CvxRad$ cannot be strictly larger than $\frac{L(A)}{\|id_V-A\|}$ for any $A$, and thus for all. 
\end{proof}

Recall that by Examples \ref{nrmex1} and \ref{nrmex2} and by Proposition \ref{grpnrmtop}, the operator norm and any composition of it with a non decreasing subadditive function is a group-norm; and that given a group-norm $N$, a left-invariant metric is obtained by $d_N(g,h)=N(g\inv h)$. With this, the group norm in the definition of a holonomic space, the usual operator norm and the convexity radius are related in the following Lipschitz condition. 
\bc\label{lpnorm} Given a holonomic space $(V,H,L)$, with $V$ an inner product space, then
$$
\|a-b\|\leq\frac{1}{\CvxRad}L(a\inv b),
$$
for all $a,b\in H$
\ec

\bt\label{holrad0} Let $(V,H,L)$ be a holonomic space and suppose further that $V$ is an inner product space and that the norm is given by $\|\cdot\|^2=\langle\cdot,\cdot\rangle$. Then the holonomy radius at the origin is given by
\be\label{holrad}
\HolRad(0)=\inf_{a\in H}\frac{L(a)}{\sqrt{2\|id_V-a\|}},
\ee
where for any $T:V\rightarrow V$,  $\|T\|$ denotes its operator norm.
\et
\begin{proof} Using the inner product, and the fact the symmetry of the group-norm $L$, $L(A\inv)=L(A)$,  and that $H$ acts by isometries, \eqref{holspeq} is equivalent to
$$
\|v-w\|^2-\|Av-w\|^2\leq L^2(A),
$$
which when expanded out yields,
$$
\|v\|^2+\|w\|^2-2\langle v,w\rangle-\|v\|^2-\|w\|^2+2\langle Av,w\rangle\leq L^2(A),
$$
and thus
$$
 2\langle Av-v,w\rangle \leq L^2(A).
$$

\noindent Thus if $\|v\|,\|w\|\leq \frac{L(A)}{\sqrt{2\|id_V-A\|}}$ then 
$$
2\langle Av-v,w\rangle \leq2\|A-id_v\|\|v\|\|w\|\leq L^2(A).
$$
Since the inequality has to hold for any $A$, it follows that
$$
\HolRad(0)\geq \inf_{a\in H}\frac{L(a)}{\sqrt{2\|id_V-a\|}}.
$$
\noindent Furthermore, for $\rho>\frac{L(A)}{\sqrt{2\|id_V-A\|}}$, let  $\varepsilon>0$ such that 
$$\varepsilon<\|id_V-A\|\rho-\frac{L^2(A)}{2\rho}=\frac{\|id_V-A\|}{\rho}\left(\rho^2-\frac{L^2(A)}{2\|id_V-A\|}\right)>0.
$$ 
\noindent By the definition of operator norm, there exists $v\in V$, with $\|v\|=\rho$ and 
$$
\|id_V-A\|\rho\geq \|Av-v\|>\|id_V-A\|\rho-\varepsilon.
$$ 
Set $w=\frac{\rho}{\|Av-v\|}(Av-v)$. It now follows that
$$
2\langle Av-v,w\rangle= 2\rho\|Av-v\|> 2\rho(\|id_V-A\|\rho-\varepsilon)> L^2(A).
$$

Thus,  \eqref{holspeq} cannot hold for $\rho>\frac{L(A)}{\sqrt{2\|id_V-A\|}}$ and the claim follows.
\end{proof}

\bc[cf. Corollary \ref{lpnorm}]\label{lipnorm} Given a holonomic space $(V,H,L)$, with $V$ an inner product space, then
$$
\sqrt{2\|a-b\|}\leq\frac{1}{\HolRad(0)}L(a\inv b),
$$
for all $a,b\in H$
\ec
\bx\label{counterex} The existence of an $r>0$ satisfying \eqref{ctid} (guaranteed for holonomic spaces by \ref{contid}) is not equivalent to the existence of an $r'>0$ satisfying \eqref{holspeq}. This follows from \eqref{holrad} by considering the following action: Let $V=\C^2$, $H=\R$, $t\cdot(z,w)=(e^{it}z,e^{\sqrt{2}it}w)$, and $L(t)=|t|$. Indeed, by Theorem \ref{holrad0}, 
$$
\inf_{a\in H}\frac{L(a)}{\sqrt{2\|id_V-a\|}}=\lim_{t\rightarrow 0^+} \frac{|t|}{\sqrt{2\sqrt{2-2\cos(\sqrt{2}t)}}}=0,
$$
whereas, by Theorem \ref{cvxrad0}, any positive $r\leq\frac{1}{\sqrt{2}}$ will make \eqref{ctid} hold. 
\ex

%%%%%%%%%%%%%%%%%%%%%%%%%%%%%%%%%%%%%%%%
\section{Holonomy}\label{Holex}
The starting point for studying the metric geometric properties of bundles over Riemannian manifolds is to consider their total spaces as Riemannian manifolds such that the projection is a Riemannian submersion. Existence and naturality of such metrics has been addressed and studied from a purely differential geometric viewpoint (see \cite{MR974641} or \cite{MR1012509} for the tangent bundle).

One procedure to view a vector bundle as  a Riemannian submersion is to endow the base with a Riemannian metric and to require that the bundle be equipped with a bundle metric and any compatible bundle connection. These two ingredients provide a plethora of metrics on the total space of the bundle (see \cite{MR2106375}), perhaps the simplest of which is the Sasaki-type metric, introduced for the tangent bundle by \citet{MR0112152}. These are just a particular case of the general construction over locally fibered maps as in \cite{MR1202431}. 

\subsection{Background on Connections and Holonomy} The results in this subsection are fairly classical (see \cite{MR1724021}).
\bd[With notation as in \cite{MR0365399}] Given a (normed) vector space $V$, there is a canonical isomorphism between $V\times V$ and $TV$, given by
\be\label{curli}
\Ii_v(w)f=\Ii(v,w)f=\dev{t}\bigg|_{t=0}f(v+tw).
\ee
That is, $\Ii_vw$ is the directional derivative at $v$ in the direction $w$. 
\ed

\br Given any vector bundle $(E,\pi)$, \eqref{curli} yields a bundle isomorphism between $\oplus^2E:=E\oplus E$ and the vertical distribution $\V=\ker\pi_*\subseteq TE$, in a natural way; that is, there is a natural transformation (also denoted by $\Ii$) from the functor $\oplus^2$ to the functor $T$.
\er

\bp\label{3map}
A connection on $(E,\pi,M)$ can be interpreted as a splitting $C$ of the following short exact sequence of bundles over the the total space $E$.
\be\label{ttses}
\xymatrix{0\ar[r] &  \pi^*E\ar[r]^{\Ii} & TE\ar[r]^{\hspace{-.4cm}\psi} & \pi^*TM\ar@/^{5mm}/[l]_{C} \ar[r]& 0}
\ee
where $\psi=(\pi_E,\pi_*)$, by regarding $C(e,u)$ as the horizontal lift of the vector $x\in M_{\pi e}$ to $e$.
\ep

Many structures on vectors bundles are transfered automatically by universality. In particular, given a vector bundle with connection $(E,\pi,C)$ over a manifold $M$, parallel translation along a curve $\alpha:I\rightarrow M$ is the trivialization of $\alpha^*E$ such that the vertical projection coincides with the projection onto the linear factor:

\bp\label{plltran} Let $(E,\pi,C)$ be a vector bundle with connection over a manifold $M$, and let $\alpha:I\rightarrow M$ be a smooth curve. The pullback becomes a bundle with connection $(\alpha^*E,\alpha^*\pi,\alpha^*C)$. Moreover, since $I$ is contractible, $\alpha^*E$ is trivial. $\alpha^*C$ yields a trivialization, called {\em parallel translation}, by considering the flow $P=P^{\alpha}$ of the vector field
\be
e\mapsto [\alpha^*C]((s,e),\dev{s})=(s,C(e,\dot\alpha)),
\ee
under the usual presentation of pullbacks as subsets of cartesian products. 
Furthermore, $P$ satisfies the following properties.
\begin{align}
P_*\del{t}&=[\alpha^*C](P,\dev{s})\\
P_t\circ P_{\tau}&=P_{t+\tau}\\
(\alpha^*\pi)\circ P(t,(s,e))&=s+t\\
P_t(e+\lambda f)&=P_te+\lambda P_tf\\
P_{t*}[\alpha^*C](e,v)&=[\alpha^*C](P_t(e),v)\\
P_{t*}\Ii(e,f)&=\Ii(P_te,P_tf)
\end{align}
so that if $I=[a,b]$, $\alpha^*E\cong I\times[\alpha^*E]_a$ by
\be
\xymatrix{
(s,e)\ar@{|->}[r]^{\mathcal{P\phantom{hhh}}}&(s,P_{a-s}e)}
\ee
\ep

\br
A metric on a vector bundle $(E,\pi)$ is a function $g:\oplus^2E\rightarrow\R$ satisfying the usual conditions. Given $(E,\pi)$ and a vector bundle with metric $(F, \tilde\pi,h)$ there is a natural metric on $\pi^*F=E\oplus F$ as a bundle over $E$ given by the pullback metric
\be
\pi^*h=h\circ(\oplus^2pr_2).
\ee
\er

\br\label{bunmet}
Given two bundles with metrics $(E, \pi_1,g)$, $(F, \pi_2,h)$ over $M$, there is a natural metric on their Whitney sum as bundles over $M$:
\be
g\oplus h=g\circ(\oplus^2pr_1)+h\circ(\oplus^2pr_2).
\ee
\qed
\er

\br\label{compconn} 
Given a vector bundle with metric and connection $(E,h,C,\nabla)$, parallel translation is by isometries.
\er

\bd Given a bundle with metric and connection, parallel translation yields a map from the space of piecewise smooth loops at a point $p\in M$, $\Omega_p$, to the group $GL(E_p)$ by
\be
\alpha\in\Omega_p\mapsto H(\alpha)=P_1^{\alpha}.
\ee 
The holonomy group $Hol_p$ at the point $p$ on the base manifold is then defined as the continuous image of $H$.
\ed

%%%%%%%%%%%%%%%%%%%%%%%%%%%%%%%%%%%%%%%
\subsection{Sasaki-type metrics} This subsection reviews the definitions of Sasaki-type metrics on general vector bundles, states a few properties and presents Theorems \ref{lengthproof} and \ref{holfib}.

In view of \ref{3map} and \ref{bunmet} there is a very natural way to define a complete Riemannian metric on the total space $E$, the Sasaki-type metric. \citet*{MR2330461} have introduced a larger class of such metrics of which the Sasaki-type is a particular case. 

\bd[\cite{MR0112152}]\label{sasmet} Given a vector bundle with metric and compatible connection $(E,\pi,h,\nabla^E)$ over a Riemannian manifold $(M,g)$, the {\em Sasaki-type metric} $\G=\G(g,h,\nabla^E)$ is defined as follows
\begin{align}
\G(e^v,f^v)&=h(e,f)\\
\G(e^v,x^h)&=0\\
\G(x^h,y^h)&=g(x,y),
\end{align}
\ed
\br An equivalent phrasing of $\G$ can be given in terms of the connection map $\kappa:TE\rightarrow E$, uniquely determined by requiring that
\be
\kappa(\sigma_*x)=\nabla^E_x\sigma;
\ee
so that $\G$ becomes
\be
\G(\xi,\eta)=g(\pi_*\xi,\pi_*\eta)+h(\kappa\xi,\kappa\eta),
\ee
for vectors $\xi,\eta\in T_eTE$.
\er

\bp \label{flatpull}Given a curve $\alpha:I\rightarrow M$ (parametrized by arc length), the (trivial) pullback bundle $\alpha^*E$ (as in \ref{plltran}) is further isometric to $I\times\R^k$ where $k$ is the rank of $E$.
\ep
\begin{proof} In view of \ref{plltran} and \ref{compconn}, one gets that
\[\alpha^*\G=\ell(\alpha)^2dt^2+\alpha^*h_p,\]
where $p=\alpha(0)$, and $\ell$ denotes the length of $\alpha$.
\end{proof}

\bp\label{distsas} The length distance on $(E,\G)$ is expressed as follows. Let $u,v\in E$
\be
d_{E}(u,v)=\inf\left\{\sqrt{\ell(\alpha)^2+\|P^{\alpha}_1u-v \|^2}\hspace{6pt} \bigg|\text{ } \alpha:[0,1]\rightarrow M,\alpha(0)=\pi u,\alpha(1)=\pi v\right\}.
\ee
Furthermore, if $\pi u=\pi v$ then
\be\label{distfor}
d_{E}(u,v)=\inf\{\sqrt{L(a)^2+\|au-v \|^2} : a\in Hol_p\},
\ee
with $L$ being the infimum of lengths of loops yielding a given holonomy element.
\ep
\begin{proof} The first expression is essentially the definition of distances in view of \ref{flatpull}. In the case $\pi u= \pi v$, the claim follows by partitioning the set of all curves $\alpha$ according to the holonomy element they generate. 
\end{proof}
\bt\label{lengthproof} Let $Hol_p$ be the holonomy group over a point $p\in M$ of a bundle with metric and connection and suppose that $M$ is Riemannian. Then the function $L_p:Hol_p\rightarrow \R$, 

\be\label{lengthnorm}
L_p(A)=\inf\{\ell(\alpha)|\alpha\in\Omega_p,P^{\alpha}_1=A\},
\ee
is a group-norm for $Hol_p$ 
\et
\begin{proof} Positivity is immediate from the fact that it is defined as an infimum of positive numbers. To prove non-degeneracy suppose that an element $A\neq I$ has zero length. There exists $u\in E_p$ such that $Au\neq u$; thus, by \eqref{distfor}, choosing $a=A$ yields $d(u,Au)=0$. A contradiction.
The length of the inverse of any holonomy element is the same because the infimum is taken essentially over the same set. Finally, to establish the triangle inequality, note that the loops that generate $AB$ contains the concatenation of loops generating $A\in Hol_p$ with loops generating $B\in Hol_p$. 
\end{proof}

\bd\label{lengthdef}
The function $L_p$, defined by \eqref{lengthnorm} will be called {\em length norm} of the holonomy group induced by the Riemannian metric at $p$.
\ed

\bt\label{holfib} Let $E_p$  be the fiber of a vector bundle with metric and connection $E$ over a Riemannian manifold $M$ at a point $p$. Let $Hol_p$ denote the associated holonomy group at $p$ and let $L_p$ be the group-norm given by \eqref{lengthnorm}. Then $(E_p,Hol_p, L_p)$ is a holonomic space. Moreover, if $E$ is endowed with the corresponding Sasaki-type metric, the associated holonomic distance coincides with the restricted metric on $E_p$ from $E$.
\et
\begin{proof}
According to the definition given in \ref{holsp}, the only remaining condition is given by \eqref{holspeq}. To see this, one needs only to note that the fiber $E_p$ is a totally geodesic submanifold of $E$. With this, given any point $u\in E_p$, let $r=\CvxRad_p(E)>0$, the convexity radius;  thus, for any pair of points $v,w\in B^E_r(p)\cap E_p$ there exists a unique geodesic from $v$ to $w$. This geodesic is necessarily $t\mapsto u-t(v-u)\in E_p$.
\end{proof}

%%%%%%%%%%%%%%%%%%%%%%%%
\subsection{Holonomy Radius of a Riemannian Manifold} Given a Riemannian manifold $(M,g)$, in view of the fundamental theorem of Riemannian Geometry, one immediately obtains a vector bundle, a connection and a bundle metric compatible with the connection; i.e. the tangent bundle, the Levi-Civita connection and the metric itself. This is the metric introduced by \citet{MR0112152}.
\subsubsection{Definitions and Basic Properties}
\bd\label{holradmfl} Let $(M,g)$ be a Riemannian manifold and let $p\in M$. The {\em holonomy radius} of $M$ at $P$ and denoted by $\HolRad_M(p)$ is defined to be the supremum of $r>0$ such that for all $u,v\in M_p$ with $\|u\|,\|v\|\leq r$ and for all $a\in Hol_p$
\be 
\|u-v\|^2-\|au-v\|^2\leq L_p^2(a),
\ee
where $L_p$ is the associated length norm on $Hol_p$.
\ed 
\br This is simply the holonomy radius at the origin of the holonomic space $(T_pM,Hol_p, L_p)$.
\er
\bt\label{posradmfl} Given a Riemannian manifold $M$. The function that assigns to each point its holonomy radius is strictly positive.
\et
\begin{proof}
This is a direct consequence of \ref{posrad} and the fact that the tangent spaces are holonomic by \ref{holfib}. 
\end{proof}
\br This fact also follows directly from geometric considerations given that $0<\CvxRad_{TM}(0_p)\leq\HolRad_M(p)$, where $\CvxRad_{TM}$ is the convexity radius of $TM$ with its Sasaki metric.
\er
\bp\label{holflat} If there exists a point $p$ in a Riemanian manifold $M$ for which the holonomy radius is not finite, then $M$ is flat. 
\ep
\begin{proof} by \ref{flatinfty}, the existence of such point is equivalent to the group being trivial. In particular, the restricted holonomy group is trivial. This is equivalent to flatness. 
\end{proof}
\br The converse is certainly not true. Consider for example a cone metric on $\R^2\setminus\{0\}$.
\er
\bc Let $M$ be a simply connected Riemannian manifold. If there is a point on $M$ with infinite holonomy radius, then $M$ is isometric to a Euclidean space.
\ec
%%%%%%%%%%%%%%%%%%%%%%%%%%%%%%%%%%
\subsubsection{Two-dimensional examples}
In the case when $(M,g)$ is a two-fold more can be said from the Gau\ss-Bonnet Theorem. Furthermore, in the particular case of the $\sS^2$ or $\hH^2$, $L$ can be computed by virtue of the isoperimetric inequality. 

Recall the following classical result. 

\bl Let $(M^2,g)$ be a $2$-dimensional Riemannian manifold and let $\gamma:[0,\ell]\subseteq\R\rightarrow M$ be any curve parametrized by arc length. Let $k$ be a signed geodesic curvature of $\gamma$ with respect to an orientation of $\gamma^*TM$.  Let $\theta(t)$ be the angle between $\dot\gamma$ and its parallel translate at time $t$. Then
\be 
2\pi-\theta(t)=\int_0^tk
\ee
Assume further that $\gamma$ is a loop. Then, possibly up to a reversal in orientation, the holonomy action of $\gamma$  at $p=\gamma(0)$ is the rotation by $2\pi-\int_0^{\ell}k$. 
\el
\begin{proof} Consider a compatible parallel almost complex structure on $\gamma^*TM$, $J$. With respect to the orthonormal frame given by $\{\dot\gamma,J(\dot\gamma)\}$, $\nabla_{\dot\gamma}\dot\gamma=kJ(\dot\gamma)$, and thus the equation for any parallel vector field $P=a\dot\gamma+bJ(\dot\gamma)$ along $\gamma$ is given by 
\begin{eqnarray*}
\dot a &=&\phantom{-}kb\\
\dot b &=&-ka
\end{eqnarray*}
which integrates to a rotation by $-\int k$ as claimed. 
\end{proof}
 \bt Let $M^2$ be a complete simply connected two-dimensional non-flat space-form with curvature $K$. Let $L:\sS^1\rightarrow\R$ be the associated length-norm on the holonomy group. Then 
 \be
L(\theta)= \dfrac{\sqrt{4\pi|\theta|\pm\theta^2}}{\sqrt{|K|} },
 \ee
 for $-\pi\leq\theta\leq\pi$, where the sign is opposite to the sign of the curvature.
 \et
 \begin{proof} By the Gau\ss-Bonnet Theorem, $ \theta=2\pi-\int k=KA$, where $A$ is the area of the region enclosed by any loop $\gamma$, so that
 $$A=\left|\frac{\theta}{K}\right|.$$ 
 Now, the isoperimetric inequality in this case  (see \cite{MR0500557}) is given by 
 $$
 \ell^2\geq4\pi A-KA^2,
 $$
 where the equality is achieved when $\gamma$ is metric circle. 
 The claim follows. 
 \end{proof}
 \bc Let $M^2$ be a simply connected two-dimensional non-flat space-form with curvature $K$. The holonomy radius at any point $p\in M$ is given by
 \be
 \inf_{-\pi\leq\theta\leq\pi}\sqrt{\dfrac{4\pi|\theta|\pm\theta^2}{2|K|\sqrt{2-2cos(\theta)}}}.
 \ee 
 \ec
 \begin{proof} In view of \eqref{holrad}, the only remain part is to compute $\|a-id\|$ for any holonomy element $a$. Since all of them are rotations by some angle $\theta$, if follows that $\|au-u\|=\|a-id\|\|u\|$ for any given $u\in T_pM$. Hence a direct application of the law of cosines yields the result. 
 \end{proof}
%%%%%%%%%%%%%%%%%%%%%%%%%%
\subsection{Length of loops---New topologies for the holonomy group}
Controlling the length of loops that generate a given holonomy element has many applications, as pointed out by \citet{MR1045885} , in Control Theory, Quantum Mechanics, or  sub-Riemannian geometry (see \cite{MR1867362}). 

Considering the infimum $L(a)$ of lengths of loops that generate a given holonomy element $a$ is a natural pick, and exhibits the fibers the vector bundle as a holonomic space as seen in \ref{holfib}. 

Although the function $a\mapsto L(a)$ is in general not even upper-semicontinuous when regarded as a function on the holonomy group with the subspace topology (or even its Lie group topology), as pointed out by \citet{MR1125864}, the following results gives a more positive outcome.

\bt\label{holtop} Let $H$ be the holonomy group of a metric connection on a vector bundle $E$ over a Riemannian manifold. There exists a finer metrizable topology on $H$ (coming from its associated length norm) so that the function $a\mapsto L(a)$ is continuous with respect to this topology and furthermore, the group action $H\times E_p\rightarrow E_p$ remains continuous. 
\et
\begin{proof} By \ref{conthol} the action map $H\times E_p\rightarrow E_p$, is continuous,  so by \ref{COcont}, the identity map is continuous from the $L$-topology to the Lie topology. Furthermore, by \ref{contnorm} $L$ is continuous with respect to the $L$-topology. 
\end{proof}

Now, the following fact hints a type of `wrong way' inheritance. 

\bp[\cite{MR1125864, MR1088509}] \label{boundedlength}Let $\pi:P\rightarrow M$ be a smooth principal bundle over a smooth manifold $M$, let a smooth connection on $\pi:P\rightarrow M$ be given, and let $H_p$ denote the holonomy group of this connection attached to some element $p$ of $P$. Suppose that $H_p$ is compact. Then there exists a constant $K$ such that every element of $H_p$ can be generated be a loop of length not exceeding $K$. 
\ep 

So, in the language of the induced length structure the following is true.
\bt\label{holbdd} Let $E\rightarrow M$ be a vector bundle with bundle metric and compatible connection. Let $H$ be the holonomy group of this connection. If $H$ is compact with the standard Lie group topology (in particular bounded with respect to any ---invariant--- metric), then $H$ with the induced length metric given by \eqref{lengthnorm} is bounded.
\et

\citet{MR1793685} introduces a way to measure the size of a holonomy transformation as a supremum over {\em acceptable} left invariant metrics. A smooth invariant metric $m$ is acceptable if for any $X\in\frak{k}=\lie(\Phi)$, the Lie algebra of $\Phi$, 
\be
\|X\|_m\leq\sup_{v,\|v\|=1} \|X(v)\|,
\ee
where $X(v)$ means the evaluation of the fundamental vector on $F$ associated with $X$. The size of a holonomy transformation $A$ is then defined as the supremum of its distances to the identity  $dist_m(A,Id)$ over acceptable metrics $m$.  And the following fact relates this `size' to the norm defined by \eqref{lengthnorm}. 
\bp[\citet{MR1793685}. Proposition 7.1]\label{Tappsize} Let $E\rightarrow B$ be a Riemannian vector bundle over a compact simply connected manifold $B$. Let $\nabla$ be a compatible metric connection and let its curvature $R$  be bounded in norm, $|R|\leq C_R$. Fix a point $x\in B$ and let $Hol(\nabla)$ be the corresponding holonomy group at $x$. Then there exists a constant $C(B)$ such that for any loop $\alpha$ in $B$, $|P_{\alpha}|\leq C\cdot C_R\cdot \ell(\alpha)$, where $P_{\alpha}\in Hol(\nabla)$ stands for the holonomy transformation induced by $\alpha$.
\ep
\bt\label{Tapphol} With the assumptions as in the previous statement, the norm given by \eqref{lengthnorm} and Tapp's holonomy size are related by $|g|\leq C\cdot C_R\cdot L(g)$, so that the induced length topology is finer than that of Tapp's holonomy size.
\et
\begin{proof} This is immediate from the inequality, since the infimum is taken over loops with the same holonomy transformation associated. 
\end{proof}
\bibliographystyle{plainnat}
\bibliography{ref}
\end{document}